\documentclass{article}  
\usepackage{latexsym, amssymb, amsmath}
\usepackage{cite}
\usepackage{mathptmx}       
\usepackage{helvet}         
\usepackage{courier}        
\usepackage{type1cm}        
\usepackage{makeidx}         
\usepackage{graphicx}        
\usepackage{multicol}        
\usepackage{subeqnarray}
\usepackage[bottom]{footmisc}
\newtheorem{definition}{Definition}[section]  
 \newtheorem{theorem}{Theorem}[section]                           
\newtheorem{remark}{Remark}[section]
\newtheorem{lemma}{Lemma}[section]

\newcommand{\h}{\hspace*{5 ex}}
\usepackage{epstopdf}
\usepackage{pdfsync, authblk}

\begin{document}

\title{Variational Inequality Approach to Stochastic Nash Equilibrium Problems with an Application to Cournot Oligopoly}
\author[1]{B. Jadamba}
\affil[1]{\em Center for Applied and Computational Mathematics, 
Rochester Institute of Technology, 85 Lomb Memorial Drive, Rochester, NY 14623, USA}
\author[2]{F. Raciti\thanks{Corresponding author. Email: fraciti@dmi.unict.it}}
\affil[2]{\em Dipartimento di
Matematica e Informatica, Universit\`a di Catania, Viale A. Doria 6-I,
95125 Catania, Italy}
\date{March 12th, 2014}
\maketitle

\begin{abstract}In this note we investigate stochastic Nash equilibrium problems by means of monotone variational inequalities in probabilistic Lebesgue spaces. We apply our approach to a class of  oligopolistic market equilibrium problems where the data are known through their probability distributions.
\end{abstract}

\section{Introduction}\label{sec_intro}

In this paper we deal with stochastic Nash equilibrum problems (SNEPs)  which we analyze using the powerful tool of stochastic variational inequalities (SVIs). As an important application of the considered SNEP model, we investigate the oligopolistic market equilibria with uncertain data.
  In the present contribution, our objective is to establish a connection between general SNEPs and SVIs and propose a model of oligopolistic markets where the cost functions are not necessarily quadratic and the demand price is not restricted  to be linear.

 We emphasize that in the deterministic framework, it is well known that oligopolistic market equilibria  are particular cases of Nash equilibria and that Nash equilibrium problems are equivalent to variational inequality problems under suitable differentiability hypotheses (see \cite{BaCa} for the infinite-dimensional case, and \cite{GaMo} for a finite-dimensional setting). Recently, some authors used a Lebesgue space formulation of oligopoly models to introduce time-dependent data (\cite{BarMau, BarDiv, BarMauro}).  We consider SNEPs where the data are affected by a certain degree of uncertainty here, for example the case that the data are known by their probabilistic measures only.  To provide a theoretical justification, we provide a formulation of SNEPs in Lebesgue spaces with probability measure, and then derive the associated SVIs.   This mechanism allows us to exploit the recently developed tools of theory of stochastic  variational inequalities in Lebesgue spaces (see e.g. \cite{GR1, GR2, GR3, GR4}). We remark that in recent years many researchers devoted their efforts to SVIs and SNEPs \cite{CheFuk:05, CheZhaFuk:09, DeXu:09, GuOzRo:99, LuBu, Pat:08, Patric, RavSha:10, RavSha:11, ShaXu:08, Xu:10}. However, these approaches differ from ours. They rely on defining a deterministic representative of the original stochastic variational inequality, and then use sample-average approximation techniques to get an estimate of the solution. In our previous work \cite{JaKhaRa:14} we have done a comparative study to analyze different approaches for a traffic equilibrum problem, and the work illustrated different solution concepts and numerical methods.

This work is organized in 7 sections.  In Section~\ref{sec2} we formulate the Nash equilibrium problem in a Lebesgue space with probability measure, and derive its equivalent stochastic variational inequality under suitable hypotheses. In Sections~\ref{sec3} and \ref{sec4} we recall some theoretical  results from \cite{GR2, GR4} together with a  description of the approximation procedure used for the solution of the stochastic variational inequality. In Section~\ref{sec5} we propose a  model of  Cournot oligopoly  with uncertain data and discuss the hypotheses needed to exploit the theory of stochastic variational inequalities. Section~\ref{sec6} is devoted to numerical examples:  we introduce  a stochastic  version of a class of utility functions widely used in the literature which yield to nonlinear monotone stochastic variational inequalities.
  The final section contains summary of results and an outline of future research directions.

\section{Stochastic Nash games and variational inequalities}\label{sec2}

Let $(\Omega,\mathcal A,P)$ be a probability space and
consider a noncooperative game with $m$ players each acting in a selfish manner in order to maximize their individual welfare. For $P-$ almost every $\omega$, each player $i$ has a strategy vector $q_{i}=(q_{i1},\ldots, q_{in})\in X_{i}(\omega)$, where $X_{i}(\omega)\subset \mathbb{R}^{n}$ is a convex and closed set,
and a utility (or welfare) function $$w_{i}:\Omega \times X_{1}(\omega) \times X_{2}(\omega)\times \ldots\times X_{m}(\omega) \to \mathbb{R}.$$ He/she choses his/her strategy vector $q_{i}$ so as to maximize $w_{i}$, given the moves $(q_{j})_{j\neq i}$ of the other players. We will use the notation 
$$q_{-i}=(q_{1},\dots,q_{i-1},q_{i+1},\ldots,q_{m}),\;\; q=(q_{i}, q_{-i}).$$

\begin{definition}  A stochastic Nash equilibrium is a random vector $q^{*}(\omega)=(q_{1}^{*}(\omega),\ldots,q_{m}^{*}(\omega))\in X(\omega)= X_{1}(\omega) \times X_{2}(\omega)\times \ldots\times X_{m}(\omega)$, such that $P-a.s.$ (almost surely):
\begin{equation}\label{parnash}
w_{i}(\omega, q_{i}^{*}(\omega),q_{-i}^{*})\geq w_{i}(\omega,q_{i},q_{-i}^{*}(\omega)), \; \forall q_{i}\in X_{i}(\omega), \forall i \in \{1,\ldots,m\}.
\end{equation}
\end{definition}

The following theorem  relates Nash equilibrium problems and variational inequalities. For its proof it suffices to apply the classical finite-dimensional proof, for each fixed value of the random parameter $\omega$.  
\begin{theorem} Let $w_{i}(\omega,\cdot)\in C^{1}(\mathbb{R}^{mn}), \forall i$, and concave with respect to $q_{i}$. Let $F: \Omega \times \mathbb{R}^{mn}\to
\mathbb{R}^{mn}$ be the mapping built with the partial gradients of the utility functions as follows:
\[
F(\omega,q)= (-\nabla_{q_{1}}w_{1}(\omega,q),\ldots, -\nabla_{q_{m}}w_{m}(\omega,q)\,).
\] 
Then, $q^{*}(\omega)\in X(\omega)$ is a stochastic Nash equilibrium if and only if, $P-a.s.$, it satisfies the variational inequality:
\begin{equation}\label{parvi}  F(\omega,q^{*}(\omega))\cdot(q-q^{*}(\omega))=
\sum_{r=1}^{m}-\nabla_{q_{r}}w_{r}(\omega,q^{*}(\omega))\cdot (q_{r}-q_{r}^{*}(\omega)) \geq 0, \forall q\in X(\omega).
\end{equation}
\end{theorem}
Problems \eqref{parnash} and \eqref{parvi} are parametric versions of the  deterministic problems, where the random parameter $\omega$ belongs to the given sample space $\Omega$. A solution $ q^{*}(\omega)$ of these problems is a random vector. 
From a statistical point of view it is important that $q^{*}(\omega)$ has  finite first and second moments. As a consequence, we introduce  integral versions of  \eqref{parnash} and \eqref{parvi}.\\
Thus, let $p\geq 2$, $\forall i $ define the set:
 $$K_{i}=\{v \in L^{p}(\Omega,P,\mathbb{R}^{n}): v(\omega) \in X_{i}(\omega), P-a.s. \}$$ and consider the problem of finding $u^{*}\in L^{p}(\Omega,P,\mathbb{R}^{mn})$ such that $
\forall i \in \{1,\ldots,m\}$ one has:
\begin{equation}\label{lebnash}
\int_{\Omega} w_{i}(\omega,u^{*}(\omega)) dP_{\omega}=
 \max_{u_{i}\in K_{i}}\int_{\Omega}w_{i}(\omega, u_{i}(\omega),u^{*}_{-i}(\omega))dP_{\omega}.
\end{equation}
The  associated variational inequality problem is the following:

Find $u^{*} \in  K_{P}$: 
\begin{equation}\label{lebvi}
\int_{\Omega}\sum_{r=1}^{m} -\nabla_{q_{r}}w_{r}(\omega, u^{*}(\omega))\cdot\,(u_{r}(\omega)-u_{r}^{*}(\omega))dP_{\omega} \geq 0\;\; \forall u\in K_{P}
\end{equation}
where
$$K_{P}= K_{1}\times \cdots \times K_{m}.$$
In \eqref{lebnash} we have introduced, $\forall i$, the functional $J_{i}: L^{p}(\Omega,P,\mathbb{R}^{mn})\to \mathbb{R}$ through:
\begin{equation}\label{lebfunc}
J_{i}(u_{1},\ldots, u_{m})=\int_{\Omega}w_{i}((\omega), u_{1}(\omega),\ldots,u_{m}(\omega))dP_{\omega}
\end{equation}
In order that this functional be well defined and to work with \eqref{lebnash} and \eqref{lebvi}, we shall impose a set of assumptions on the functions $w_{i}$, namely:
\begin{enumerate}
\item[(a)] For all $i\in\{1,\ldots,m\}$, $w_{i}(\cdot,q)$ be a random variable with respect to the sigma-algebra defined on $\Omega$, $\forall q$, and $w_{i}(\omega, \cdot)\in C^{1}(\mathbb{R}^{mn})$ P-a.s.
\item[(b)] For all $i\in\{1,\ldots,m\}$, $w_{i}(\omega,0)\in L^{1}(\Omega, P)$.
\item[(c)] For all $i\in\{1,\ldots,m\}$, $w_{i}(\omega, q)$ be concave with respect to $q_{i}$, $P-a.s.$, for all fixed values of $q_{-i}$.
\item[(d)] For all $i\in\{1,\ldots,m\}$, $|\nabla_{q}w_{i}(\omega,q)| \leq \alpha(\omega) + \beta_{i}(\omega) |q|^{p-1}$,
where $\beta_{i} \in L^{\infty}(\Omega,P)$ and $\alpha \in L^{p'}(\Omega,P)$, $p'=p/p-1$.\label{growth}
\end{enumerate}
We are now in a position to prove a simple lemma which is fundamental for the sequel.
\begin{theorem} Let assumptions a)-d) be fulfilled. Then, for all $i\in\{1,\ldots,m\}$ the functional $J_{i}(u)=J_{i}(u_{i},u_{-i})$ is well defined on $L^{p}(\Omega,P, \mathbb{R}^{mn})$ and concave with respect to the variable $u_{i}\in L^{p}(\Omega,P, \mathbb{R}^{n})$ for each fixed $u_{-i}$. Moreover $J_{i}$ is Gateaux-differentiable with respect to $u_{i}$, for each $u_{-i}$ and its derivative is given by:
\begin{equation}\label{gat}
D_{i}J_{i}(u_{i},u_{-i})(v_{i})= \int_{\Omega} \nabla_{q_{i}}w_{i}(\omega, u_{i}, u_{-i})\, \cdot v_{i}(\omega) dP_{\omega}=
\end{equation}
\begin{equation}{\nonumber}
\int_{\Omega}\sum_{r=1}^{n} \left( \frac{\partial}{\partial q_{r}} \,(w_{i}(\omega, u(\omega))\right)v_{ir}(\omega) dP_{\omega},\;\; \forall v_{i}=(v_{i1},\ldots,v_{in}) \in L^{p}(\Omega,P, \mathbb{R}^{n})
\end{equation}
\end{theorem}
{\bf Proof}. \ 
First, we  show that the functional $J_{i}$ is well defined for all $i$. Thus, for P-almost every $\omega \in \Omega$ apply
Lagrange Theorem to the function $w_{i}(\omega,q)$, with respect to the interval of endpoints $0,q$. We get that $\exists \xi \in \mathbb{R}^{mn}$, $|\xi| < |q|$ such that:
$$| w_{i}(\omega,q)| \leq |w_{i}(\omega,0)| + |\nabla_{q}w_{i}(\omega,\xi)| |q|\leq
|w_{i}(\omega,0)| +|\alpha (\omega)| |q| + \beta_{i}(\omega)|\xi|^{p-1} |q|.$$
Then, $\forall u\in L^{p}(\Omega, P, \mathbb{R}^{mn})$ we get
$$ |w_{i}(\omega, u(\omega))| \leq |w_{i}(\omega,0)| + |\alpha(\omega)| |u(\omega)|^{p-1} +\beta_{i}(\omega)| u(\omega)|$$
which shows that $w_{i}(\omega,u(\omega))$ belongs to $L^{1}(\Omega,P)$, hence $J_{i}$ is well defined. The concavity of $J_{i}(u_{i},u_{-i})$ with respect to $u_{i}$ is a straightforward consequence of the analogous property  di
$w_{i}(\omega,q)$.\\
In order to prove that $J_{i}$ is Gateaux-differentiable with respect to $u_{i}$, for each fixed $u_{-i}$ , fix a point $u_{i}$, a direction $v_{i}$ and for each $t\in ]0,1[$ consider the quotient:
\begin{eqnarray*}
&&\frac{J_{i}(u_{i}+tv_{i}, u_{-i})-J_{i}(u_{i},u_{-i})}{t}\\
&=& \int_{\Omega} \frac{1}{t}\left[  w_{i}(\omega, u_{i}(\omega)+tv_{i}(\omega),u_{-i}(\omega))- w_{i}(\omega, u_{i}(\omega), u_{-i}(\omega))\right]  dP_{\omega}\\
&=&\int_{\Omega} \nabla_{q_{i}}w_{i}(\omega, u_{i}(\omega) +th(\omega)v_{i}(\omega), u_{-i}(\omega))\cdot v_{i}(\omega)    dP_{\omega}, 
\end{eqnarray*}
where $h: \Omega \to [0,1]$ is  a random variable. We obtain \eqref{gat} because it is possible to pass to the limit under the integral sign for $t\to 0$. Indeed, since $w_{i}(\omega, \cdot)$ has continuous partial derivatives, it follows that for $t\to 0$,  we get,  $P-a.s.$:
$$\nabla_{q_{i}} w_{i}(\omega, u_{i}(\omega) +th(\omega)v_{i}(\omega), u_{-i}(\omega))\cdot v_{i}(\omega)
\longrightarrow     \nabla_{q_{i}} w_{i}(\omega, u_{i}(\omega), u_{-i}(\omega))\cdot v_{i}(\omega),$$ 
moreover
\begin{eqnarray*}
&&|\nabla_{q_{i}} w_{i}(\omega, u_{i}(\omega) +th(\omega)v_{i}(\omega), u_{-i}(\omega))\cdot v_{i}(\omega)|\\
&\leq&
  |\alpha (\omega)| |v_{i}(\omega)| + \beta_{i}(\omega)\,\left(\;|u_{i}(\omega)|+|v_{i}(\omega)| + |u_{-i}(\omega)|\right)^{p-1}.
\end{eqnarray*}
At last,  the fact that $D_{i}J_{i}(u)(\cdot)$ is  a linear and continuous functional on $L^{p}(\Omega,P, \mathbb{R}^{n})$
concludes the proof. 
\hfill $\square$\\
Once we have established the expression of the Gateaux derivative of $J_{i}$, consider, for each $u$, the operator $\Gamma (u):L^{p}(\Omega, P, \mathbb{R}^{mn})\to L^{p'}(\Omega, P, \mathbb{R}^{mn})$ defined by:
$$\Gamma(u)=(-D_{1}J_{1}(u), \ldots, -D_{m}J_{m}(u)). $$
Then, from the infinite dimensional theory of Nash equilibrium problems,  we get (see e.g. \cite{BaCa}) that \eqref{lebnash} is equivalent to
$$ u^{*}\in K_{P}: \Gamma(u^{*})(u-u^{*}) \geq 0,\; \forall u \in K_{P},$$
which is nothing other than \eqref{lebvi}.
\section{Stochastic variational inequalities in Lebesgue spaces}\label{sec3}
In the sequel we shall study SNEPs and, in particular, the oligopolistic market, through its equivalent variational inequality \eqref{lebvi}. As mentioned in the introduction, variational inequalities of this kind have been introduced  quite recently and in this section we recall  the main results useful for our application.  A more comprehensive treatement can be found in \cite{GR2, GR3, GR4}. In particular, we shall treat the case where the deterministic and random variable are separated and in this case an approximation procedure for the computation of the solution is presented.  

Let $(\Omega,{\mathcal A},P)$ be a probability space. Let $G,H:\mathbb{R}^{k} \rightarrow \mathbb{R}^{k}$ be two given maps, let  $b,c \in \mathbb{R}^k$ be fixed vectors, and let $R$ and $S$ be two real-valued random variables defined on $\Omega$. Let $\lambda$ be a random vector in $\mathbb{R}^k,$ let $D$ be random vector in $\mathbb{R}^m,$ and let $A\in
\mathbb{R}^{m\times k}$ be a given matrix. For $\omega \in \Omega,$ we define a random set
 $$M(\omega):= \{x \in \mathbb{R}^k :\  Ax \leq D(\omega) \}.$$
Consider the following stochastic variational inequality: For almost all $\omega \in \Omega,$ find $\hat{x}:=\hat{x}(\omega)\in M(\omega)$
 such that
 \begin{equation} \label{sepa} \langle S(\omega)\,  G(\hat{x}) + H(\hat{x}),
 z- \hat{x}\rangle \geq
  \langle R(\omega)\, c+b, z- \hat{x}\rangle,\quad \textrm{for every}\ z \in M(\omega).
 \end{equation}
Variational inequality \eqref{sepa} holds pointwise on $\Omega,$ except a fixed null set depending on the solution $\hat{x}$. To facilitate the foregoing discussion, we set
  $$F(\omega,x) := S(\omega)\,G(x) + H(x),$$
Let $S,G$ and $H$ be such that $F :  \Omega \times  {\mathbb{R}^k}
\mapsto \mathbb{R}^k$ is a  Carath\'eodory function. That is, for each fixed $x \in \mathbb{R}^k$, the function
$F(\cdot,x)$ is measurable with respect to $\mathcal A$ whereas for each $\omega \in \Omega\,$ the function $F (\omega,\cdot)$ is  continuous. We also assume that $F(\omega,\cdot)$ is monotone for every $\omega\in \Omega$:
$$\langle F(\omega,x)-F(\omega, y),x-y \rangle \geq 0, \; \forall x,y, \forall \omega$$
 If the equality sign  holds only for $x=y$, then $F$ is said strictly monotone and, in this case there is at most a solution to \eqref{sepa} which, under
    suitable conditions  belongs to an $L^p$ space for some $p\geq 2.$ 
    
    A stronger form of monotonicity will be useful in the sequel:
    \begin{definition} F is strongly monotone, uniformly with respect to $\omega$ iff $\exists a >0$:
   $$\langle F(\omega,x)-F(\omega, y),x-y \rangle \geq  a\|x-y\|^{2}, \; \forall x,y, \forall \omega.$$
    \end{definition}
For this, we proceed to derive the integral formulation of \eqref{sepa}. For a fixed $p\geq 2$, we define the reflexive Banach  space $L^{p} (\Omega,P, \mathbb{R}^k)$ of random vectors $V$ from $\Omega$ to $ \mathbb{R}^k$ such that the expectation ($p$-moment) is given by:
$$ E^{P}\|V\|^{p} =\int_{\Omega} \|V(\omega)\|^{p} dP(\omega) <\infty.$$
For the subsequent development, we need the following growth condition
\begin{equation}\label{gro}\|F(\omega,z)\|
\leq \alpha(\omega) +\beta(\omega) \|z\|^{p-1},\quad\forall z\in \mathbb{R}^k, \quad \textrm{for some}\ p\geq 2,
\end{equation}
where $\alpha\in L^{p'}(\Omega,P)$ and $\beta\in L^{\infty}(\Omega,P).$

Due to the above growth condition, the Nemitsky operator $\hat{F}$ associated to $F$, acts from $L^p(\Omega, P,\mathbb{R}^k)$ to  $L^{p'}(\Omega, P,\mathbb{R}^k),$ where $p^{-1}+{p'}^{-1}=1.$  Furthermore, we have $$\hat{F}(V)(\omega):=F(\omega,V(\omega)),\quad \omega\in \Omega.$$
Assuming $D \in L_m ^{p} (\Omega):= L^{p} (\Omega, P,
\mathbb{R}^m),$ we introduce the following nonempty, closed and convex subset of $L_k ^{p} (\Omega)$
 $$M^{P}:= \{  V \in  L_k^{p} (\Omega) : A \,
 V(\omega)  \leq D(\omega),\; P- a.s. \},$$
 which is the $L^p$ analogue of $M(\omega)$ defined above.

Let $S(\omega)\in L^{\infty},$  $0< \underline s<S(\omega) < \overline s,$ and $R(\omega)\in L^{p'}$. Equipped with these notation, we consider the following $L^p$ formulation of \eqref{sepa}.  Find $\hat U \in M^P$  such that for every $V\in M^P,$ we have
\begin{equation} \label{intsepa}
\int_{\Omega} \langle S(\omega) \, G(\hat U (\omega)) +
H(\hat U(\omega)), V(\omega) - \hat U(\omega) \rangle \, dP(\omega) \geq
\int_{\Omega} \langle b+ R(\omega) \, c, V(\omega) - \hat U(\omega)\rangle dP(\omega).
\end{equation}
To get rid of the abstract sample space $\Omega$, we consider the joint distribution $\mathbb{P}$ of the random vector $(R,S,D)$ and work with the special probability space
 $(\mathbb{R}^{d}, {\mathcal B} (\mathbb{R}^{d}), \mathbb{P})$,
 where the dimension $d:= 2+m$. For simplicity, we assume that $R$, $S$ and $D$ are independent random
 vectors. We set
 \begin{eqnarray*}
 r&=&R(\omega),\\
 s&=&S(\omega),\\
 t&=&D(\omega),\\
 y &=& (r,s,t).
 \end{eqnarray*}
 For each
$y \in \mathbb{R}^d$, we define the set
  $$ M(y) :=  \{ x \in \mathbb{R}^k :\  Ax \leq t \}.$$
The pointwise formulation of the variational inequality reads: Find $\hat x $ such that $\hat x(y) \in M(y)$,
$\mathbb{P}$ - a.s., and
the following inequality holds for $\mathbb{P}$ - almost every
$y \in \mathbb{R}^{d}$ and for every $x\in M(y),$ we have
 \begin{equation}\label{prob}
\langle s\, G(\hat x (y))  +  H(\hat x (y)),  x-\hat x (y) \rangle \geq
\langle  r c+b, x-\hat x (y) \rangle \,.
\end{equation}

In order to obtain the integral formulation of \eqref{prob}, consider the space $L^{p} (\mathbb{R}^{d}, \mathbb{P},\mathbb{R}^k)$ and introduce the closed and convex set
 $$M_{\mathbb{P}}:= \{ v \in L^{p} (\mathbb{R}^{d},
        \mathbb{P},\mathbb{R}^k): A v (r,s,t) \leq t, \;\mathbb{P}-a.s. \}.$$
With this terminology, we consider the variational inequality of finding  $\hat u \in M_{\mathbb{P}}$ such that for every $v \in
M_{\mathbb{P}}$ we have
 \begin{equation} \label{intprob}
\int\limits_{0}^{\infty} \int\limits_{\underline s}^{\overline s}\int\limits_{\mathbb{R}^d} \langle s \, G(\hat u (y))  + H(\hat u (y)),
 v(y) - \hat u(y)\rangle  \, d\mathbb{P}(y) \geq
\int\limits_{0}^{\infty} \int\limits_{\underline s}^{\overline s}\int\limits_{\mathbb{R}^d} \langle b + r \, c, v(y)-\hat u(y)\rangle\, d\mathbb{P}(y).
\end{equation}
\begin{remark}
Our approach and analysis extends readily to more general finite Karhunen-Lo\`{e}ve expansions
$$
\lambda (\omega)=b + \sum_{l=1}^L R_l(\omega)~c_l,\quad
F(\omega,x) = H(x) + \sum_{l=1}^{L_F} S_l(\omega)\,G_l(x).$$
\end{remark}
%
\section{An Approximation Procedure by Discretization of Distributions}\label{sec4}
Without any loss of generality, we assume that $R \in L^q(\Omega, P)$ and $D \in L_m^{p} (\Omega, P)$  are nonnegative (otherwise we can use the standard decomposition
 in the positive part and the negative part).
Moreover, we assume that the support, the set of possible outcomes, of $S \in L^{\infty} (\Omega , P)$ is the interval
$[{\underline s}, \overline s ) \subset (0, \infty )$. Furthermore, we assume that the probability measures $P_R$, $P_S$, and $P_D$ are
continuous with respect to the Lebesgue measure, so that according to
the theorem of Radon-Nikodym, they have the probability densities $\varphi _R$, $\varphi _S$, and $\varphi _{D_i},$ $i=1, \dots, m $, respectively.
Therefore,  for  $i=1,\dots, m,$ we have
\begin{eqnarray*}
\mathbb{P} &=&P_R \otimes P_S \otimes P_D,\\
dP_R (r)&=& \varphi _R (r)\, dr,\\
dP_S (s)&=& \varphi _S (s)\, ds\\
dP_{D_i}(t_i) &=& \varphi_{D_i}(t_i)\,dt_i.
\end{eqnarray*}
Notice that $v \in L^{p} (\mathbb{R}^d, \mathbb{P},
\mathbb{R}^k )$ means that $(r,s,t) \mapsto \varphi_R (r)  \varphi_S (s) \varphi_D (t) v (r,s,t)$
belongs to the Lebesgue space
$L^{p} (\mathbb{R}^d,\mathbb{R}^k)$ with respect to the Lebesgue measure where
$$\varphi_D (t) := \prod_i  \varphi_{D_i}(t_i).$$
Therefore, we can define the probabilistic integral variational inequality: Find $\hat u:=\hat{u}(y) \in M_{\mathbb{P}}$ such that for every
$v \in M_{\mathbb{P}}$, we have
 \begin{eqnarray*}
\int\limits_{0}^{\infty} \int\limits_{\underline s}^{\overline s}
 \int\limits_{\mathbb{R}_{+}^m}
 \langle s \, G(\hat u) + H(\hat u),  v- \hat u \rangle \,
        \varphi_R (r)  \varphi_S (s) \varphi_D (t) \, dy \geq
        \int\limits_{0}^{\infty} \int\limits_{\underline s}^{\overline s}
        \int\limits_{\mathbb{R}_{+}^m}
\langle b + r \,  c, v-\hat u \rangle \, \varphi_R (r)  \varphi_S
        (s) \varphi_D (t) \, dy \,.
\end{eqnarray*}
For numerical approximation of  the solution $\hat u$, we begin with a discretization of the space $X:= L^{p} (\mathbb{R}^d,\mathbb{P},\mathbb{R}^k).$ For this, we introduce a sequence  $\{ \pi _n \}_n$ of
partitions of the support
$$\Upsilon := [0, \infty) \times [\underbar   s,\overline
s) \times \mathbb{R}_{+}^m $$ of the probability measure $\mathbb{P}$
 induced by the random elements $R,S,$ and $D$.
For this, we set
$$\pi _n = (\pi_n ^R , \pi_n ^S,
\pi_n ^D),$$
where
  \begin{eqnarray*}
        \pi_n^R &:=& (r_n^0, \dots, r_n^{N_n^R}),\\
        \pi_n^S &:=& (s_n^0, \dots, s_n^{N_n^S}),\\
        \pi_n^{D_i} &:=& (t_{n,i}^0, \dots, t_{n,i}^{N_n^{D_i}})\\
        0& =&r_n^0 <r_n^1 < \dots r_n^{N_n^R}=  n\\
        \underbar s &=&s_n^0 <s_n^1 < \dots s_n^{N_n ^S}= \overline s\\
0 &=&t_{n,i}^0 <t_{n,i}^1 < \dots t_{n,i}^{N_n^{D_i}}=  n \, \;
        (i= 1, \dots, m)\\
        |\pi _n ^{R} |& :=& \max \{ r_n^j-r_n^{j-1}: j=1,\dots, N_n^{R} \}
        \rightarrow 0 \h (n \rightarrow \infty )\\|\pi _n^{S} | &:=&\max \{ s_n^k-s_n^{k-1}: k=1,\dots, N_n^S  \}
        \rightarrow 0 \h (n \rightarrow \infty )\\
        |\pi _n ^{D_i} | &:=& \max \{ t_{n,i}^{h_i}-t_{n,i}^{h_i -1}:
        h_i=1,\dots, N_n^{D_i} \} \rightarrow 0 \h (i= 1, \dots, m;\,
         n \rightarrow \infty) \,.
\end{eqnarray*}
These partitions give rise to the exhausting sequence $\{\Upsilon_n\}$
of subsets of $\Upsilon$, where each $\Upsilon_n$ is given by the finite
disjoint union of the intervals:
  $$
        I_{jkh}^n :=  [r_n^{j-1}, r_n^{j})
         \times [s_n^{k-1}, s_n^{k}) \times I_h^n \,,
  $$
where we use the multi-index $h = (h_1,\cdots,h_m)$ and
  $$I_h ^n := \displaystyle \Pi_{i=1}^m  \,
\displaystyle [t_{n,i}^{h_i-1},t_{n,i}^{h_i}).$$
For each $n \in \mathbb{N},$ we consider the space of the
$\mathbb{R}^l$-valued step functions ($l \in \mathbb{N}$)
on $\Upsilon_n$, extended by $0$ outside of $\Upsilon_n$:
  $$
        X_n^l   := \{v_n:  v_n (r,s,t)= \sum _j \sum _k \sum _h
         v^n_{jkh} 1_{I^n_{jkh}} (r,s,t) \,, v^n_{jkh} \in \mathbb{R}^l \}
  $$
where $1_I$ denotes the  $\{ 0,1\}$-valued characteristic
function of a subset $I$.

To approximate an arbitrary function $w \in L^{p} (\mathbb{R}^d, \mathbb{P},
\mathbb{R}),$ we employ the mean value truncation operator $\mu _0 ^ n$
associated to the partition $\pi _n $ given by
 \begin{equation}
        \mu _0^n w :=  \sum_{j=1}^{N_n^{R}}\sum_{k=1}^{N_n^{S}}
         \sum _{h}(\mu _{jkh} ^n w)\,1_{I_{jkh}^n}\,,
 \end{equation}
where
\begin{eqnarray*}
 \mu _{jkh} ^n w  :=  \left\{
 \begin{array}{ll}
 \displaystyle \frac{1}{\mathbb{P}(I_{jkh})}
        \int_{I_{jkh}^n} w(y) \, d \mathbb{P} (y)
   & \mbox{if } \mathbb{P}(I^n_{jkh})  > 0 \, ;
 \\[1ex]
  0 & \mbox{otherwise.}
\end{array} \right.
\end{eqnarray*}
Analogously, for a $L^{p}$ vector function $v=(v_1,\dots,v_l)$, we define
$$\mu _0^n
v := (\mu _0^n v_1, \dots, \mu _0^n v_l).$$ 
The basic property of the mean value truncation operator is expressed in the following lemma (see \cite{GR2}).
\begin{lemma} \label{medie}
For any fixed $l \in \mathbb{N}$, the linear operator $\mu_0^n : L^{p}
(\mathbb{R}^d, \mathbb{P}, \mathbb{R}^l) \rightarrow L^{p} (\mathbb{R}^d,
\mathbb{P}, \mathbb{R}^l)$ is  bounded with $\| \mu _0^n\| =1$ and for $n
\rightarrow \infty$, $\mu _0^n$ converges pointwise in $L^{p} (\mathbb{R}^d,
\mathbb{P}, \mathbb{R}^l)$ to the identity.
\end{lemma}
To construct approximations for
  $$M_{\mathbb{P}}= \{ v \in L^{p} (\mathbb{R}^d,
        \mathbb{P},\mathbb{R}^k): A v (r,s,t) \leq t\,,
        \;\mathbb{P}-a.s.\},$$
we introduce the orthogonal projector
$q: (r,s,t) \in \mathbb{R}^d \mapsto t  \in \mathbb{R}^m$
and define for each elementary cell $I_{jkh}^n$,
$$
{\overline q}_{jkh}^n  = ( \mu _{jkh}^{n} q) \in \mathbb{R}^m ,\;\;
         \h(\mu _{0}^{n} q) = \sum _{jkh} {\overline q}_{jkh}^n \, 1_{I^n_{jkh}} \in
        X_n^m \,.
$$
This leads to the following sequence of convex and closed sets of the polyhedral type:
  $$M_{\mathbb{P}}^n :=  \{v \in X_n^k : A v_{jkh}^n \leq
        {\overline q}_{jkh}^n \,, \; \forall j,k,h \}.$$
It has been proven (see \cite{GR2}) that the sequence
$\{ M_{\mathbb{P}}^n\}$ approximate the set $ M_\mathbb{P}$  in the sense of Mosco (\cite{Mos:69}).
In order to to approximate the random variables $R$ and $S,$ we introduce
\begin{eqnarray*}\rho_n& =&\sum_{j=1}^{N_n^R} r_n ^{j-1}\, 1_{[r_n^{j-1}, r_n ^j )} \in  X_n\\
\h\sigma_n &=&\sum_{k=1}^{N_n^S} s_n ^{k-1}\, 1_{[s_n^{k-1}, s_n ^k )}
              \in  X_n,
\end{eqnarray*}
where
\begin{eqnarray*}
\sigma _n (r,s,t) &\rightarrow& \sigma (r,s,t)=s,\quad \textrm{in}\ \ L^{\infty} (\mathbb{R}^d, \mathbb{P})\\
\rho_n (r,s,t) &\rightarrow&
\rho(r,s,t)=r,\quad \textrm{in}\ \ L^{p} (\mathbb{R}^d, \mathbb{P}).
\end{eqnarray*}
Combining the above ingredients, for $n \in \mathbb{N}$, we consider the following discretized variational inequality: Find $\hat u _n:=\hat{u}_n(y) \in M_{\mathbb{P}}^n $
such that for every $v_n \in
M_{\mathbb{P}}^n$, we have
  \begin{equation}
 \label{Pn}  \int\limits_{0}^{\infty} \int\limits_{\underline s}^{\overline s}\int\limits_{\mathbb{R}^d}
 \langle \sigma_n(y) \, G(\hat u_n)  + H(\hat u_n),
v_n - \hat u_n  \rangle \, d\mathbb{P}(y) \geq
   \int\limits_{0}^{\infty} \int\limits_{\underline s}^{\overline s}\int\limits _{\mathbb{R}^d} \langle b + \rho_n(y) \, c,
  v_n - \hat u_n\rangle\, d\mathbb{P}(y) \,.
\end{equation}
It turns out that \eqref{Pn} can be split in a finite number of finite dimensional variational inequalities: For every $n \in
\mathbb{N},$ and for every $j,k,h$ find $  \hat u^n_{jkh} \in M^n_{jkh} $ such
that
  \begin{equation}\label{split}
        \langle\tilde F_{k}^n (\hat u^n_{jkh}),
        v^n_{jkh}- \hat u^n_{jkh}\rangle \geq
       \langle \tilde c^{n}_{j}, v^n_{jkh}- \hat u^n_{jkh}\rangle, \quad \textrm{for every}\ v^n_{jkh} \in M^n_{jkh},
  \end{equation}
where
\begin{eqnarray*}
        M^n_{jkh}  &:=&  \{v^n_{jkh} \in \mathbb{R}^k :
        A v^n_{jkh} \leq {\overline q}_{jkh}^n \} \,, \\[1ex]
        \tilde F_{k}^n &:=&
        s_n^{k-1} \, G + H\\
         \tilde c_{j}^n \,&:=& \, b +  r_n^{j-1} \, c.
  \end{eqnarray*}
Clearly, we have
$$ \hat u_n =  \sum _j \sum _k \sum _h \hat u^n_{jkh} \, 1_{I^n_{jkh}}
         \in X_n^k.$$
Now, we can state the following convergence result
(whose proof can be found in \cite{GR4}).
\begin{theorem} \label{CONV} Assume that $F(\omega,\cdot)$ is strongly monotone uniformly with respect to $\omega\in \Omega.$ Then
 the sequence $ \hat u_n $  generated by the substitute problems in  \eqref {Pn}
converges strongly in $L^{p} (\mathbb{R}^d, \mathbb{P},\mathbb{R}^k)$
for $n \rightarrow \infty$ to the unique solution $\hat u$ of \eqref {intprob}.
\end{theorem}
\begin{remark} Looking carefully at the proof in \cite{GR4}, we can deduce that if the uniform strong monotonicity hypothesis is not satisfied, and  $F$ is only monotone, but we know that the solution is unique we obtain weak convergence of $ \hat u_n $ to $\hat u$. This implies convergence of the approximate mean values to the exact mean value of the solution.
\end{remark}
\section{The stochastic oligopoly model}\label{sec5}
In this section, we propose a model of  oligopolistic market with uncertain data and show that the theoretical and numerical tools developed in the previous sections can be successfully applied to the model.  
The classical oligopolistic market equilibrium problem is a Nash game with a special structure and it was first  introduced by A. Cournot \cite{Co} a long time ago.  Recent  years  have witnessed a renewed interest in oligopoly theory, and  many specific  cases of oligopolistic markets have been studied in detail, for instance the electricity market (see e.g.  \cite{ChHoal, ChSial}, and   \cite {BoOgAlMa} for  a model based on real industrial data). 

We consider here the case in which $m$ players are the producers  of the same commodity. The quantity produced by firm $i$ is denoted by $q_{i}$ so that $q\in \mathbb {R}^{m}$ denotes the global production vector.  Let $(\Omega,\mathcal A, P)$ be a probability space and for every  $i\in\{ 1,\ldots m \}$ consider functions $f_{i}:\Omega \times \mathbb{R} \to \mathbb{R}$ and  $p:\Omega \times \mathbb{R}^{m} \to \mathbb{R}$.
More precisely,  $f_{i}(\omega, q_{i})$ represents the cost of producing the commodity for firm $i$, and is assumed to be, $P-a.s.$, nonnegative, increasing and $C^{1}$, while   
$p(\omega, q_{1}+\ldots+q_{m})$ represents the demand price associated with the commodity. For $P-$almost every $\omega \in \Omega$,  $p$ is assumed nonnegative, increasing and $C^{1}$.  The resulting welfare function $w_{i}$ is assumed to be concave with respect to $q_{i}$.
We also assume  that all these functions are random variables w.r.t. $\omega$, i.e. they are measurable with respect to the probability measure $P$ on $\Omega$. In this way, we cover the possibility that both the production cost and the demand price are affected by a certain degree of uncertainty, or randomness. 
Thus, the welfare (or utility) function of player $i$, representing the net revenue,  is given by:
\begin{equation} 
w_{i}(\omega, q_{1},\ldots,q_{m})= p(\omega,q_{1}+\ldots+q_{m})q_{i} -f_{i}(\omega, q_{i}).
\end{equation}

Although many models assume no bounds on the production, in a more  realistic model the production capability is bounded from above and we also allow these upper bounds to be  random variables: $0\leq q_{i}\leq {\overline q}_{i}(\omega)$. 
Thus, the specific Nash equilibrium problem associated with this model takes the following form:\\
 For $P-a.e.$ $\omega\in \Omega$, find $q^{*}(\omega)=\left(q_{1}^{*}(\omega),\ldots,q_{m}^{*}(\omega)\right)$:
\begin{equation}\label{parcour}
w_{i}(\omega,q^{*}(\omega))=\max_{0\leq q_{i}\leq {\overline q}_{i}(\omega)} 
\{p(\omega, q_{i}+ \sum_{j\neq i} q_{j}^{*}(\omega)\,)q_{i}, \;   -f_{i}(\omega,q_{i})\},\; \forall i\in\{1,\ldots,m\}. 
\end{equation}
In order to write the equivalent variational inequality, consider,$\forall \omega$, a  closed and convex subset of $\mathbb{R}^{m}$:
$$
K(\omega)= \{(q_{1},\ldots, q_{m}): 0\leq q_{i} \leq {\overline q}_{i}(\omega),\,\forall i \}
$$
and define the  functions
\begin{eqnarray} 
F_{i}(\omega,q)&:=&  \frac{\partial f_{i}(\omega,q_{i})}{\partial q_{i}} -
 \frac{\partial p(\omega,\sum_{j=1}^{m} q_{j})}{\partial q_{i}}q_{i}- p(\omega,\sum_{j=1}^{m} q_{j})\\
&=& f'_{i}(\omega,q_{i})- p'(\omega,Q)q_{i}-p(\omega,Q),\;\;\; (Q=\sum_{j=1}^{m} q_{j}). \nonumber
\end{eqnarray}
The Nash problem is then equivalent to the following variational inequality:
 for $P-$a.e. $\omega\in \Omega$, find $q^{*}(\omega) \in K(\omega)$ such that
\begin{equation}\label{parvicour}
\sum_{i=1}^{m} \left[ \frac{\partial f_{i}(\omega,q_{i}^{*}(\omega))}{\partial q_{i}} -
 \frac{\partial p(\omega,\sum_{j=1}^{m} q_{j}^{*}(\omega))}{\partial q_{i}}q_{i}- p(\omega,\sum_{j=1}^{m} q_{j}^{*}(\omega))\right](q_{i}-q_{i}^{*}(\omega)) \geq 0
\end{equation}
$\forall q\in K(\omega)$. 
 \begin{remark}
 Since $F(\omega,\cdot)$ is continuous, and $K(\omega)$ is convex and compact, problem \eqref{parvicour} is solvable for almost every $\omega \in \Omega$, due to the Stampacchia's theorem. In the case that the production capability is assumed unbounded some additional hypotheses (i.e. coercivity, see  e.g. \cite{MauRac}) have to be present to ensure the solvability of \eqref{parvicour}. 
  \end{remark}
   Moreover, we assume that $F(\omega, \cdot)$ is monotone, i.e.:
$$
\sum_{i=1}^{m} (F_{i}(\omega,q)-F_{i}(\omega, q') )(q_{i}-q_{i}\rq{})\geq 0\ \; \forall \omega \in \Omega,\;\; \forall q, q' \in \mathbb{R}^{m}.  
$$
(recall that $F$ is said to be strictly monotone if the equality holds only for $q=q'$ and in this case \eqref{parvicour} has a unique solution).
It is noteworthy that some classes of utility functions widely used in the economic literature enjoy some form of monotonicity (see section~\ref{sec6}). 

Now we are interested in computing statistical quantities associated with the solution $q^{*}(\omega)$, in particular its mean value. For this purpose, in accordance with the general scheme of Section 2,  we consider a Lebesgue space formulation of problems  \eqref{parvicour}: Find $u^{*}\in K$ such that
\begin{eqnarray}{\nonumber}\int_{\Omega}\label{lebvicour}\sum_{i=1}^{m} \left[ \frac{\partial f_{i}(\omega,u_{i}^{*}(\omega))}{\partial q_{i}} -
 \frac{\partial p(\omega,\sum_{j=1}^{m} u_{j}^{*}(\omega))}{\partial q_{i}}u_{i}- p(\omega,\sum_{j=1}^{m} u_{j}^{*}(\omega))\right]\times \\
(u_{i}(\omega)-u_{i}^{*}(\omega)\,)dP_{\omega}\geq 0\;\; ,
\end{eqnarray}
where
$$K=\{ u\in L^{p}(\Omega,P, \mathbb{R}^{m}): 0\leq u_{i}(\omega)\leq {\overline q_{i}}(\omega) \},\;\overline q_{i}\in L^{p}(\omega,P).   $$
Since the stochastic  oligopolistic market problem will be studied through \eqref{lebvicour}, we ensure its solvability by the following theorem:
\begin{theorem}\label{teovinashleb}
Let $f_{i}(\cdot,q_{i}), p(\cdot,\sum_{j=1}^{m}q_{j})$ be measurable, and $f_{i}(\omega,\cdot), d_{i}(\omega, \cdot)$ be of class $C^{1}$. Let $F$ be  strictly monotone and satisfy the growth condition $d)$ of Section~\ref{sec2}. Then \eqref{lebvicour} admits a unique solution.
\end{theorem}
{\bf Proof}. \ 
Under our assumptions, $F:\Omega \times \mathbb{R}^{m} \to \mathbb{R}^{m}$ is a Carath\'eodory function and it is well known that for each measurable function $u(\omega)$, the function $F(\omega, u(\omega))$ is also measurable. Under the growth condition $d)$ the superposition operator $\hat{F}: u(\omega) \to F(\omega, u(\omega))$ maps $L^{p}(\Omega,P,\mathbb{R}^{m})$ in $L^{p'}(\Omega,P,\mathbb{R}^{m})$ and is continuous, being $P$ a probability measure. Moreover, the uniform strong monotonicity of $F$ implies the strong monotonicity of $\hat{F}$. The set $K_{P}$ is convex, closed and (norm) bounded, hence weakly compact. Then, monotone operator theory applies  and  \eqref{lebvicour} admits a unique solution (see e.g. \cite{MauRac} for a recent survey on existence theorems which also includes the case of unbounded sets).
\hfill $\square$\\
Now, in view of the numerical approximation of the solution, we further  specialize our model and  assume that the random and the deterministic part of the operator can be separated. Thus, we assume that the price can be affected by two  random perturbations $\alpha (\omega)$ and $S(\omega)$ such that:
 $$p(\omega, Q)=S(\omega)p(Q)+\alpha(\omega),$$
 while the cost functions are of the type:
 $$f_{i}(\omega,q_{i})=\beta_{i}(\omega)f_{i}(q_{i})+g_{i}(q_{i}),$$
 that is, the cost functions consists of a deterministic term $g_{i}$ and a term, (still denoted by $f_{i}$ with an abuse of notation), which is {\em modulated} by the random perturbation $\beta_{i}$.  
 Here $\alpha, \beta_{i} $ are real random variables, with  
 $0<\underline s\leq S(\omega)\leq {\overline s}\;$,
 $0<\underline \beta_{i}\leq \beta_{i}(\omega)\leq {\overline \beta_{i}}$. 

As a consequence, the operator $F$ takes the form:
$$
F_{i}(\omega,q)=\beta_{i}(\omega)\frac{\partial f_{i}(q_{i})}{\partial q_{i}}+  \frac{\partial g_{i}(q_{i})}{\partial q_{i}}  -S(\omega)p(\sum_{j=1}^{m}q_{j})-\alpha(\omega)-S(\omega)\frac{\partial p(\sum_{j=1}^{m}q_{j})}{\partial q_{i}}q_{i}.
$$
 Furthermore, we assume that $F$ is uniformly strongly monotone, $w_{i}(\omega,0) \in L^{1}(\Omega)$,  and  the growth condition $d)$ of Section~\ref{sec2} is satisfied.

   Now, according to the methodology explained in  Section~\ref{sec3},  we will   work with the probability distributions induced on the images of the functions: $A = \alpha(\omega),s=S(\omega), B_{i} =\beta_{i}(\omega), Q_{i}={\overline q}_{i}(\omega)$. Thus, let $y=(A,s,B,Q)$ and consider the probability space $(\mathbb {R}^{d}, \mathcal{B}, \mathbb{P})$ with $d=2+2m$, where ${\mathcal B}$ is the Borel sigma-algebra on  $\mathbb {R}^{d}.$
In order to formulate our problem in the image space, we introduce the closed convex set $K_{\mathbb P}$ by: 
$$\;K_{\mathbb P}=\{  u\in L^{2}(\mathbb{R}^{d},\mathbb{P}, \mathbb{R}^{m}): 0\leq u_{i}(A,s,B,Q) \leq \overline{Q_{i}}, \forall i, \mathbb{P}-a.s.\}\;.$$
We assume that all the random variables are independent and that each probability distribution is characterized by its density $\varphi$. Thus, we have 
$\mathbb{P}= \mathbb{P}_{A}\otimes \mathbb{P}_{s}\otimes \mathbb{P}_{B}\otimes  \mathbb{P}_{Q}$, $dP_{\alpha }(A)=\varphi_{\alpha}(A)dA,\;dP_{S }(s)=\varphi_{S}(s)ds,\;
 dP_{\beta}(B)=\varphi_{\beta}(B)dB,\;
dP_{\overline q}(Q)=\varphi_{\overline q}(Q)dQ$,  where we used the compact notation 
$ \displaystyle\varphi_{x}(X)= \prod_{i=1}^{n} \varphi_{x_{i}}(X_{i})$. 
Thus, we obtain the following problem: Find $u^{*}\in K_{\mathbb P}$ such that $\forall u \in K_{\mathbb P}$
\begin{eqnarray}
\nonumber
&&\int\limits_{\underline s}^{\overline s}\int\limits_{\underline \beta}^{\overline \beta} \int\limits_{\mathbb{R_{+}}}\int\limits_{\mathbb R^{m}_{+}}\sum_{i=1}^{m}
 \Bigg[B_{i}\frac{\partial f_{i}(u_{i}^{*}(A,s,B,Q))}{\partial q_{i}}+  \frac{\partial g_{i}(u_{i}^{*}(A,s,B,Q))}{\partial q_{i}}  -p\left(\sum_{j=1}^{m}u_{j}^{*}(A,s,B,Q)\right)\\
\label{densvinash}
\nonumber
&&-A-\frac{\partial p(\sum_{j=1}^{m}u_{j}^{*}(A,s,B,Q))}{\partial q_{i}} u_{i}^{*}(A,s,B,Q)\Bigg]\times\\
&& (u_{i}(A,B,Q)- u_{i}^{*}(A,s,B,Q))\varphi_{\alpha}(A)\varphi_{S}(s) \varphi_{\beta}(B) \varphi_{\overline q}(Q)\; ds\,dA\,dB\,dQ\geq 0,
\end{eqnarray}
where the symbol $\int\limits_{\underline \beta}^{\overline \beta}$ represents the $m$ integrals $\int\limits_{\underline \beta_{i}}^{\overline \beta_{i}}\;$. This  formulation is suitable for the approximation procedure based on  discretization and truncation explained in Section~\ref{sec4}. 
\section{A class of utility functions and numerical examples}\label{sec6}
In this section we consider a modified and random version of a class of utility functions introduced by Murphy, Sheraly and
Soyster in  \cite{MuSheSo}  and successively used by other scholars.  These functions generate a nonlinear monotone variational inequality on a certain $L^{p}$ space, where $p$ is determined by the power law of the cost functions. 
 The cost and demand price functions for the five-firm case in \cite{MuSheSo} are given by:
\begin{eqnarray*}f_{i}(q_{i})&=&c_{i}q_{i}+\frac{b_{i}}{b_{i}+1} k_{i}^{-1/b_{i}}q_{i}^{\frac{b_{i}+1}{b_{i}}},\; i=1,\ldots,5\\
p(Q)&=& 5000^{1/1.1}Q^{-1/1.1}, \;\; Q=\sum_{i=1}^{5}q_{i}.
\end{eqnarray*}
The values of the parameters $c_i, k_i,b_i$ in \cite{MuSheSo}  alongwith our upper bounds for the $q_i$ are given  Table~\ref{params_example1}. An approximate solution of the problem obtained by a projection method is  given in \cite{Nag} as  $(q_1, q_2, q_3, q_4, q_5) = (36.937, 41.817, 43.706, 42.659, 39.179)$.
\\

\begin{table}[h!]
\begin{center}
\caption{Parameter values for the numerical example}
\label{params_example1}
\begin{tabular}{lccccc}
\hline
$i$ & 1 & 2 & 3 & 4& 5 \\
\hline
$c_i$ &10 & 8 & 6 & 4 & 2  \\
$k_i$ & 5 & 5 & 5 & 5 & 5  \\
$b_i$ & $1.2$ & $1.1$ & $1.0$ & $0.9$ & $0.8$  \\
$\overline{q}_{i}$&100&100&100&100&100\\
\hline
\end{tabular}
\end{center}
\end{table}

Before introducing random parameters in the above functions we note that  the demand price becomes unbounded when the total quantity $Q$ approaches $0$ (commodity is scarce). Although the solution $Q^{*}=0$ is never met in most examples, in order to deal with a well behaved function we  consider the functional  form:
$$p(Q)= 5000^{1/1.1}(Q+e)^{-1/1.1},$$
where $e$ is a small positive parameter which determines the maximum price the consumer can pay when the commodity is very scarce. In our model, we add a random perturbation $r(\omega)$ to $c_{i}$, and
we modulate the price function by a random function $S(\omega)$.

Thus,  for the general case of $m$ firms, we introduce cost functions given by:
\begin{equation}\label{modicost}
f_{i}(\omega,q_{i})= [c_{i} + r(\omega)]q_{i} + \frac{b_{i}}{b_{i}+1} k_{i}^{-1/b_{i}}q_{i}^{\frac{b_{i}+1}{b_{i}}},
 \end{equation}
 where $b_{i}, c_{i},k_{i}$ are positive parameters, and demand price functions:
 \begin{equation}\label{modiprice}
 p(\omega,Q)= \left[S(\omega)\right]^{a}\, \frac{1}{(Q+e)^{a}},
 \end{equation}
 where $0<\underline{s}<S(\omega) <\overline{s}$, and  $a$ is a parameter such that  $0<a<1$ ( $a=1/1.1$ in \cite{MuSheSo}).
 
With these functions we can build the Carath\'eodory function $F$ which defines the  variational inequality  through:
\begin{equation}\label{operand}
F_{i}(\omega,q)= c_{i} + r(\omega)  +k_{i}^{-1/b_{i}}q_{i}^{1/b_{i}}+ a [S(\omega)]^{a}  \frac{q_{i}}{{(Q+e)}^{a+1}}
-\frac{[S(\omega)]^{a}} {{(Q+e)}^{a}},\; i=1\ldots m.
\end{equation}
We also use the notation $F_{i}(\omega,q)= G_{i}(\omega,q)+H_{i}(\omega,q)$, where $G_{i}$ represents the sum of the first three terms in \eqref{operand}, while $H_{i}$ is the rest of the sum, which contains the price function. 
The monotonicity of $F$ is analyzed in the following: 
\begin{theorem} The function $F(\omega,\cdot)$ defined by \eqref{operand} is strictly monotone in $\mathbb{R}^{m}_{+}$,  for all $\omega \in \Omega$ and  for all fixed values of the parameters therein.
\end{theorem}
{\bf Proof}. \  
Let us observe that the functions $k_{i}^{-1/b_{i}}q_{i}^{1/b_{i}}$ are strictly increasing for all $i$, hence the operator $G(\omega,\cdot)$ is strictly monotone on $\mathbb{R}^{n}_{+}$, for all $\omega$.\\
In order to study the monotonicity properties of $H$, we preliminary notice that the function $p:\Omega \times \mathbb{R}\to \mathbb{R}$ defined in \eqref{modiprice} has strictly positive second derivative (w.r.t. $Q$):
$$p''(\omega,Q)=\frac{a\, (a+1)[S(\omega)]^{a} }{(Q+e)^{a+2},}$$
 (recall that $0<\underline{s}\leq S(\omega)$), therefore $p(\omega, \cdot)$ is strictly convex for all $Q\geq 0$.\\
 Let us now consider the function $\displaystyle Q\,p(\omega,Q)=[S(\omega)]^{a}\frac{Q}{(Q+e)^{a}}$ which is strictly concave on 
  $\mathbb{R}^{n}_{+}$, for each value of $\omega$, with second derivative given by:
  $$[Qp(\omega,Q)]''=   [S(\omega)]^{a}\frac{a\,(a-1)Q-2ae}{(Q+e)^{a+2}}<0$$
  ($0<a<1$). Hence, we get:
  \begin{equation}\label{conca}
  -2p'(\omega,Q)> Qp''(\omega,Q)
  \end{equation}
  which will be exploited in the sequel.\\
  To prove the strict monotonicity of $H(\omega,\cdot)$ for all $\omega$ we compute its Jacobian matrix:
\begin{eqnarray*}
  J_{ij}(\omega,Q)&=&-p'(\omega,Q)-q_{i}p''(\omega,Q),\ \textrm{if}\ i\neq j\\
  J_{ii}(\omega,Q)&=&-2p'(\omega,Q)-q_{i}p''(\omega,Q).
\end{eqnarray*}  
It is useful to decompose $J$ as follows:
  $$J(\omega,Q)=-p'(\omega,Q){\bf 1}-p'(\omega,Q)I-p''(\omega,Q)(q_{i})_{ij},$$
  where $\bf 1$ denotes the $m\times m$ matrix with each entry equal to $1$, $I$ is the identity matrix and the matrix
  $(q_{i})_{ij}$ has each entry of the row $i$ equal to $q_{i}$.
  We prove that $J(\omega,Q)$ is positive definite for all $\omega$ and for all $Q\geq 0$ by studying  the quadratic form
  $$T(\omega,Q)(h)=h^{T}J(\omega,Q)h, \; h\in \mathbb{R}^{m}.$$
  From the decomposition of the $J(\omega,Q)$ we then get:
  \begin{eqnarray*}
 T(\omega, Q)(h)&=&-p'(\omega,Q)\left(  \sum_{i,j=1}^{m} ({\bf 1})_{ij}h_{i}h_{j}+ \sum_{i,j=1}^{m}(I)_{ij}h_{i}h_{j}  \right)-
  p''(\omega,Q)\sum_{i,j=1}^{m}(q_{i})_{ij}h_{i}h_{j}\\
 &=& -\left\{p'(\omega,Q)\left[ \left(\sum_{j=1}^{m}h_{j}\right)^{2}+\|h\|^{2}\right] +    p''(\omega,Q) 
   \left(\sum_{j=1}^{m}h_{j}\right)\left(\sum_{j=1}^{m}q_{j}h_{j}\right)  \right\}.
  \end{eqnarray*}
  Now, if $h\neq (0,\ldots,0)$, from \eqref{conca} we get the strict inequality:
  $$2T(\omega,Q)(h)>p''(\omega,Q)\left\{Q\left[ \|h\|^{2}+ \left( \sum_{j=1}^{m}h_{j}\right)^{2} \right]-
  2   \left(\sum_{j=1}^{m}h_{j} \right)\left(\sum_{j=1}^{m}q_{j}h_{j} \right) \right\}.$$
  Given that $p(\omega,\cdot)$ has strictly positive second derivative, it suffices to prove that the quantity in curly brackets is nonnegative.\\
 Thus, let $h\in \mathbb{R}^{m}$ with $\sum_{j=1}^{m}h_{j}\,\geq 0$ (the case  where $\sum_{j=1}^{m}h_{j}\,\leq 0$ can be analyzed along the same lines), so that:
 $$\left(\sum_{j=1}^{m}h_{j}\right) \left(\sum_{j=1}^{m}q_{j}h_{j}\right)\leq  h_{j_{max}}
 \left(\sum_{j=1}^{m}h_{j}\right) \left(\sum_{j=1}^{m}q_{j}\right)\leq
 Qh_{j_{max}}  \left(\sum_{j=1}^{m}h_{j}\right),  $$
 where $h_{j_{max}}$ is and index such that $h_{_{max}}\geq h_{j},\; \forall j=1\ldots m$ and without loss of generality we can assume $j_{max}=m$ in the sequel. We get then:
\begin{eqnarray*}
2T(\omega,Q)(h) &>&Q\,p''(\omega,Q) \left\{ \left[ \|h\|^{2}+ \left( \sum_{j=1}^{m}h_{j}\right)^{2} \right]- 2h_{m} \sum_{j=1}^{m}h_{j} \right\}\\
&=& Qp''(\omega,Q) \left[\sum_{j=1}^{m-1}h_{j}^{2}+h_{m}^{2}+  \left(\sum_{j=1}^{m-1}h_{j}\,+ h_{m} \right)^{2}
  -2h_{m}\sum_{j=1}^{m}h_{j} \right]\\
&=&Qp''(\omega,Q) \left[\sum_{j=1}^{m-1}h_{j}^{2} + \left(\sum_{j=1}^{m-1}h_{j} \right)^{2} \right] \geq 0.
\end{eqnarray*}
Thus, $T(\omega,Q)(h) >0, \forall \omega \in \Omega, \forall Q\in \mathbb{R}^{m}_{+}, \,\forall h\neq (0,\ldots,0)$.
\hfill $\square$\\
Now, let us consider the case $m=5$ with  the data as in Table~\ref{params_example1}. The function $F$, defines a Nemitsky operator between Lebesgue spaces, as explained in the previous sections. To be precise, since the exponents $b_{i}$ in the cost functions vary from $0.8$ to $1.2$, we select $p= 1+1/0.8$ so that the Nemitsky operator associated to $F$ maps functions  $u\in L^{9/4}$ into $u\in L^{9/5}$. 
Moreover, we let random parameters $r(\omega)$ and $S(\omega)$ to have truncated normal distributions as follows:
\begin{eqnarray*}
r &\sim& -0.5 \leq N(0, 0.25)\leq 0.5\\
s & \sim& 4950\leq N(5000, 10) \leq 5050
\end{eqnarray*}
while fixing parameter $e$ at $0.0001$. Mean values
 $E(u)$ of $u(r,s)=(u_{1}, u_{2}, u_{3}, u_4, u_{5})$
  obtained by numerical  approximations are presented in Table~\ref{results_example1} where $n_r$ and $n_s$ stand for number of discretization points for intervals $[-0.5, 0.5]$ and $[4950, 5050]$ respectively. 
\begin{table}[h!]
\caption{Mean values of $u_i$, i=1,\ldots,5}
\label{results_example1}
\begin{tabular}{lccccc}
\hline\noalign{\smallskip}
 & $E(u _1)$ & $E(u _2)$ & $E(u_3)$ & $E(u_4)$& $E(u_5)$ \\
\noalign{\smallskip}\hline\noalign{\smallskip}
$(n_r,n_s)=(200, 20000)$ &36.8855 &     41.7615&      43.6448&      42.5972&       39.121 \\
$(n_r,n_s)=(400, 40000)$ & 36.913&      41.7928&      43.6776&      42.6294&      39.1506 \\
\noalign{\smallskip}\hline
\end{tabular}
\end{table}
\section{Conclusions}\label{sec7}
In this article we considered Nash equilibrium problems in Lebesgue spaces with probability measure and derived their equivalent variational inequality formulation. As a specific application, we proposed a model of oligopolistic market with uncertain data
to which the recent theory of  random variational inequality (\cite {GR2, GR4}) was applied. We also illustrated our model and the approximation procedure by means of a class of utility functions which yield to nonlinear monotone random variational inequalities.

 Further developments of our approach can be done in several directions: other type of probabilistic constraints could be considered instead of the ``robust''  pointwise constraints (see e.g. \cite{DenRus:03});  an extension of our numerical method, for example through parallelization, is desirable and would permit the treatement of problems with a larger number of independent random variables; at last, the  theory and computation of the stochastic Lagrange multipliers associated to SNEPs  in Lebesgue spaces is a topic that has been adressed only recently (\cite {JaRa}) in a simplified model and only from a theoretical point of view.


\begin{thebibliography}{40}

\bibitem{BaCa} C.~Baiocchi,  A.~Capelo,  Variational and Quasivariational Inequalities,  Applications to Free Boundary Problems, John Wiley and Sons (1984)

\bibitem{BarMau} A. Barbagallo, A. Maugeri, Duality theory for the dynamic oligopolistic market equilibrium problem, Optimization, 60, 1-2, pp. 29-52 (2011)

\bibitem{BarDiv} A. Barbagallo, R. Di Vincenzo, Lipschitz continuity and duality for dynamic oligopolistic market equilibrium problem with memory term, Journal of Mathematical Analysis and Applications, 382, 231-247 (2011)

\bibitem{BarMauro} A. Barbagallo, P. Mauro, Evolutionary variational formulation for oligopolistic market equilibrium problems with production excesses, Journal of Optimization Theory and Applications 155, 288-314 (2012)

\bibitem{BoOgAlMa} F.~Bonenti, G.~Oggioni, E.~Allevi, G.~Marangoni,  Evaluating the EU ETS impacts on
profits, investments and prices of the Italian electricity market,  Energy
Policy, 59, 242-256 (2013)

\bibitem{CheFuk:05} X. Chen, M. Fukushima, Expected residual minimization method
for stochastic linear complementarity problems, Math. Oper. Res., 30, 1022-1038 (2005)

\bibitem{CheZhaFuk:09} X. Chen, C. Zhang, M. Fukushima, Robust solution of monotone stochastic linear complementarity problems, Math. Prog. B, 117, 51-80 (2009)

\bibitem{ChHoal} Y. Chen, B. F. Hobbs, S. Leyffer, T. S. Munson,  Leader-follower equilibria for electric power and NO$_x$  allowences markets, Computational Management Science, 3(4), 307-330 (2006)  

\bibitem{ChSial} Y. Chen, J. Sijm, B. F. Hobbs, W. Lise, Implications of CO$_2$ emissions trading for short-run electricity market outcomes in northwest Europe,  Journal of Regulatory Economics, 34, 23-44 (2008)

\bibitem{Co} A. A. Cournot, Researches into the Mathematical Principles of the Theory of Wealth, 1838, English Translation, MacMillan, London, England (1897)

\bibitem{DeXu:09} V. De Miguel, H. Xu, A stochastic multiple-leader Stackelberg model:
Analysis, computation, and application, Oper. Res., 57, 1220-1235 (2009)

\bibitem{DenRus:03} D. Dentcheva, A. Ruszczy\'{n}ski, Optimization  with stochastic dominance
constraints, SIAM J. Optim., 14, 548-566 (2003)

\bibitem{FalRac:07}  F. Raciti, P. Falsaperla, Improved non-iterative algorithm
for the calculation of the equilibrium in the traffic network problem, J. Optim. Theory. Appl., 133, 401-411 (2007)


\bibitem{GaMo} D. Gabay, H. Moulin, On the uniqueness and stability of Nash Equilibria in non cooperative games, in {Applied Stochastic Control in Econometrics and Management Sciences}, editors: A.~Bensoussan, P.~Kleindorfer, and C.S.~Tapiero, 271-294, North Holland, Amsterdam (1980)

\bibitem{GuOzRo:99}  G. G\"urkan, A.Y. \"Ozge, S. M. Robinson, {Sample-path solution of stochastic variational inequalities}, Math. Program., 84,  313-333 (1999)

\bibitem{GR1} J. Gwinner, F. Raciti,  {Random Equilibrium problems on
networks},   Mathematical and Computer Modelling, 43, 880-891 (2006)

\bibitem{GR2} J. Gwinner, F. Raciti, {On a Class of Random
Variational Inequalities on Random Sets}, {Numerical Functional Analysis and
Optimization}, 27, 5-6, 619-636 (2006)

\bibitem{GR3} J. Gwinner, F. Raciti,   {On Monotone  Variational
Inequalities with  Random Data}, Journal of Mathematical Inequalities, 3 (3),  443-453 (2009)

\bibitem{GR4} J. Gwinner, F. Raciti, {Some equilibrium problems under uncertainty and random variational inequalities}, {Annals of Operations Research}, 200, 299-319 (2012),  
 DOI 10.1007/s10479-012-1109-2

 \bibitem{JaRa} B. Jadamba, F. Raciti, {On the modelling of some environmental games with uncertain data},  Journal of Optimimization Theory and Applications, DOI 10.1007/s10957-013-0389-2

\bibitem{JaKhaRa:14} B.~Jadamba, A.~A.~Khan, F.~Raciti, Regularization of  Stochastic Variational Inequalities and a Comparison of an $L_p$ and a Sample-Path Approach,  Nonlinear Analysis A,  94, 65-83 (2014), http://dx.doi.org/10.1016/j.na.2013.08.009


\bibitem{LuBu} S. Lu and A. Budhiraja. Confidence regions for stochastic variational inequalities, 
Mathematics of Operations Research, 38 (3), 2013, http://dx.doi.org/10.1287/moor.1120.0579


 \bibitem{MuSheSo} F. H. Murphy, H. Sheraly, A.L. Soyster, A mathematical programming approach for determining oligopolistic market equilibrium, Mathematical Programming, 24, 92-106 (1982)


 \bibitem{MauRac} A. Maugeri, F. Raciti, On Existence Theorems for Monotone and Nonmonotone Variational Inequalities, Journal of Convex Analysis, 16, 3\&4, 899-911 (2009)

\bibitem{Mos:69} U. Mosco, Convergence of convex sets and of solutions of
variational inequalities, Advan. Math.,  3, 510-585 (1969)
 

\bibitem{Nag} A. Nagurney, Network Economics: A Variational
Inequality Approach, Second and Revised Edition, Kluwer
Academic Publishers, Dordrecht, The Netherlands (1999)

\bibitem{Patric} M. Patriche, Equilibrium of Bayesian fuzzy economies and
quasi-variational inequalities with random
fuzzy mappings,  Journal of Inequalities and Applications, 2013:374 (2013),
doi:10.1186/1029-242X-2013-374

\bibitem{Pat:08} M. Patriksson, On the applicability and solution of bilevel optimization models in
transportation science: A study on the existence, stability and
computation of optimal solutions to stochastic mathematical
programs with equilibrium constraints, Trans. res., B 42, 843-860 (2008)

\bibitem{ShaDenRus:09}   A. Shapiro, D. Dentcheva, A. Ruszczy\'{n}ski,
Lectures on Stochastic Programming - Modeling and Theory,  SIAM, Philadelphia (2009)

\bibitem{ShaXu:08}
A, Shapiro, H. Xu, Stochastic mathematical programs with equilibrium constraints,
modelling and sample average approximation,  Optimization, 57, 395-418 (2008)

 \bibitem{RavSha:10} U. Ravat, U. V. Shanbhag, On the characterization of solution sets of smooth and nonsmooth stochastic Nash games, Proceedings of the American Control Conference (ACC), Baltimore (2010)

\bibitem{RavSha:11} U. Ravat,  U. V. Shanbhag,
 On the characterization of solution sets of
smooth and nonsmooth convex
stochastic Nash games, SIAM Journal of Optimization, 21 (3), 1168-1199 (2011)

\bibitem{RaSha} U. Ravat, U.V. Shanbhag, On the existence of solutions to
stochastic variational inequality and complementarity problems, 
arXiv:1306.0586v1 [math.OC] (2013)

\bibitem{Xu:10}
H. Xu, Sample average approximation methods for a class of
stochastic variational inequality problems, 
Asia-Pacific J. Oper. Res., 27, 103-119 (2010)




%








%


%







\end{thebibliography}
\end{document}